\documentclass[a4paper,11pt,twoside]{article}
\usepackage{latexsym}
\usepackage{amsmath}
\usepackage{amsfonts}
\usepackage{amssymb}
\usepackage{vmargin}
\usepackage{url}
\usepackage{tikz}
\usepackage{graphicx}
\usepackage{color}
\usepackage[round]{natbib}
\usepackage{times}
\usepackage[hang,flushmargin]{footmisc}
\setmarginsrb{2.5cm}{2.5cm}{2.5cm}{2.5cm}{0mm}{10mm}{0mm}{10mm}

\usepackage{fancyhdr}
\usepackage{lastpage}
\fancypagestyle{title}{%
  \setlength{\headheight}{0pt}%
  \fancyhf{}
  \fancyfoot[R]{\footnotesize \copyright The authors 2011. Published by Oxford University Press on behalf of the Institute of Mathematics and its Applications.\\
   All rights reserved. For permissions, please email: journals.permission@oup.com}
  \fancyhead[L]{\footnotesize \textsl{Teaching Mathematics and Its Applications} (2011) \textbf{30}, \thepage--\pageref{LastPage}\\
                doi:10.1093/teamat/hrr008 \quad Advance Access publication \ 12 May 2011}
  \fancyhead[R]{}
}%

\begin{document}
\title{\vspace*{-1cm} \bfseries Mollweide's formula in teaching trigonometry}
\author{\scshape Natanael Karjanto{\footnotesize$^\ast$}\\
\textsl{\small Department of Mathematics, University College, Sungkyunkwan University, Natural Science Campus}\\
\textsl{\small 2066 Seobu-ro, Jangan-gu, Suwon, 16419, Gyeonggi-do, Republic of Korea} \\ [5pt]
{\footnotesize $^\ast$\textsl{Email: \url{natanael@skku.edu}}} \\[5pt]
\rm {\footnotesize [Submitted December 2010; accepted March 2011]}}
\date{}
\maketitle
\thispagestyle{title}

\begin{abstract}
\noindent
Trigonometry is one of the topics in mathematics that the students in both high school and pre-undergraduate levels need to learn. Generally, the topic covers trigonometric functions, trigonometric equations, trigonometric identities and solving oblique triangles using the Laws of Sines and Cosines. However, when solving the oblique triangles, Mollweide's formula is most likely to be omitted from the discussion. Mollweide's formula---which exhibits a cyclical nature---is particularly useful in checking one's result after solving an oblique triangle since all six components of the triangle are involved. It is interesting to note that proving Mollweide's formula can be performed without words. Furthermore, the Law of Tangents can be derived straightforwardly from this equation. In this article, we revisit Mollweide's formula and provide classroom examples where this equation comes into alive. In addition, we suggest that this seemingly less-known equation is to be included in the mathematics syllabus on the topic of Trigonometry. \bigskip \\
{\bf Keywords:} Mollweide's formula; Trigonometry; oblique triangle; the Law of Sines; the Law of Cosines.
\end{abstract}

\section{Introduction}

The Department of Applied Mathematics at the University of Nottingham Malaysia Campus serves as a service department in teaching mathematics modules for the Faculty of Engineering students, both at the Foundation (pre-undergraduate) and Undergraduate levels. One of the topics that we would like to focus in this article is Trigonometry. Previously, Trigonometry is taught at the Foundation level in HG1FND Foundation Mathematics module. Starting from the academic year of 2009/2010, Trigonometry would be delivered under the new module F40FNA Foundation Algebra.

The materials covered in Trigonometry amongst others are trigonometric ratios, solving right triangles, trigonometric functions, fundamental trigonometric identities, sketching trigonometric functions, sum and difference trigonometric identities, double-angle and half-angle identities, product-to-sum and sum-to-product identities, trigonometric equations and additional topics on solving oblique triangles. Two basic tools for solving oblique triangles are the Law of Sines and the Law of Cosines. The derivation of these laws is readily available in many textbooks on Trigonometry, such as \cite{Barnett06} and \cite{Larson97}. The Law of Sines is used to solve triangles given two angles and any side (\textsc{asa} or \textsc{aas}) or two sides and an angle opposite one of them (\textsc{ssa}). The Law of Cosines is used to solve triangles when two sides and the included angle (\textsc{sas}) or three sides (\textsc{sss}) are given. 

The other formulations which are basically rather untouched are the Law of Tangents and Mollweide's formula. Mollweide's formula, sometimes also referred to as Mollweide's  equation, is a set of two relationships between sides and angles in a triangle~\citep{Sullivan88}. This equation is particularly useful in checking one's result after solving an oblique triangle since all six components of the triangle are involved~\citep{Wilczynski14}. Interestingly, from Mollweide's formula, we can derive the Law of Tangents almost straightforwardly. Mollweide's formula possesses symmetry property and---similar to the Laws of Sines and Cosines---exhibits a cyclical nature. This means, to go from one form of the equation to the other, we simply rotate the order of the sides and the angles, each time replacing a letter by the next in line~\citep{Holt07}.

The purpose of this article is to revisit Mollweide's formula and to propose the topic to be included in teaching Trigonometry. This article is organized as follows. In the following section, we will discuss Mollweide's formula, including the derivation of this equation, its cyclical form and the relevant examples in solving oblique triangles. Finally, the final section concludes this article and gives remarks to our discussion.

\section{Mollweide's formula}

Mollweide's formulas are given as follows:
\begin{eqnarray}
  \frac{\sin \frac{1}{2}(\alpha - \beta)}{\cos \frac{1}{2}\gamma} &=& \frac{a - b}{c} \label{1}\\
  \frac{\cos \frac{1}{2}(\alpha - \beta)}{\sin \frac{1}{2}\gamma} &=& \frac{a + b}{c}. \label{2}
\end{eqnarray}
Throughout this article, the following convention is adopted. In a $\triangle \; ABC$, $\alpha = \measuredangle \; BAC$, $\beta = \measuredangle \; ABC$ and $\gamma = \measuredangle \; ACB$, thus $\alpha + \beta + \gamma = \pi$. Furthermore, $a = BC$, $b = AC$ and $c = AB$.

The equations adopt their name from a German mathematician and astronomer Karl Brandan Mollweide (1774--1825).\footnote[7]{\noindent Mollweide is more well known in the field of cartography due to his contribution in map projection, called the Mollweide projection. More information on Mollweide projection can be found in a thorough and well-organized book on the history of map projection~\citep{Snyder97}.} Nonetheless, this pair of equations was discovered earlier by Isaac Newton (1643--1727) and fully developed by Thomas Simpson (1710--1761). An excellent overview of the history of Mollweide's formula is given by~\cite{Wu07}. Some authors have attempted to prove this equation without words. See for instance~\cite{Kleine88} and \cite{Wu01,Wu02}. The readers who are interested in the general topic on proof without words are encouraged to consult~\cite{Nelsen93,Nelsen00}.
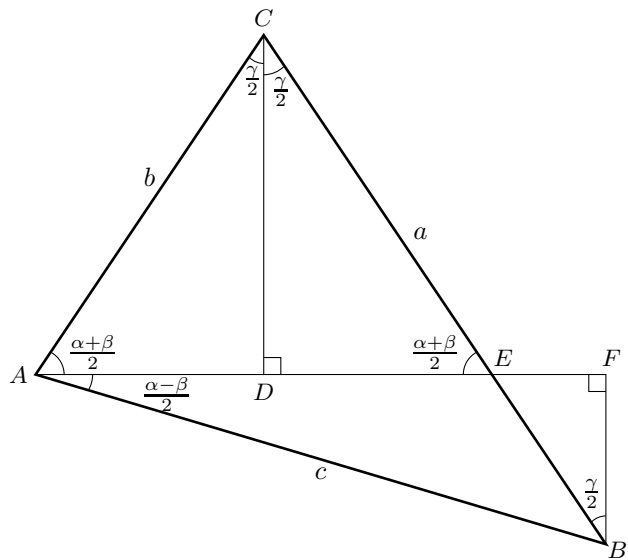
\begin{center}
\begin{tikzpicture}[scale = 0.75	]
  \draw[line width=1pt]  (0,0)--(10,-3)--(4,6)--cycle;
  \draw (0,0)--(10,0)--(10,-3);
  \draw (9.7,0)--(9.7,-0.3)--(10,-0.3);
  \draw (4,6)--(4,0);
  \draw (4,0.3)--(4.3,0.3)--(4.3,0);
  
  \draw (0.5,0) arc (0:60:4.5mm);
  \draw (1,0) arc (0:-30:5.5mm);
  \draw (10,-2.5) arc (90:135:3.5mm);
  \draw (7.5,0) arc (180:120:4.5mm);
  \draw (4,5.5) arc (-90:-135:3.5mm);
  \draw (4,5.3) arc (-90:-45:5mm);
  
  \path(-0.3,0) node(A)    {\footnotesize $A$};
  \path(10.2,-3.1) node(B) {\footnotesize $B$};
  \path(4,6.3) node(C) 	   {\footnotesize $C$};
  \path(4,-0.3) node(D)    {\footnotesize $D$};
  \path(8.2,0.3) node(E)   {\footnotesize $E$};
  \path(10.1,0.3) node(F)  {\footnotesize $F$};

  \path(6.75,2.5) node(a) {\small $a$};
  \path(2,3.5) node(b)    {\small $b$};
  \path(5,-1.75) node(c)  {\small $c$}; 
  
  \path(1,0.4) node()     {\small $\frac{\alpha + \beta}{2}$};
  \path(2.3,-0.35) node() {\small $\frac{\alpha - \beta}{2}$};
  \path(7,0.4) node() 	  {\small $\frac{\alpha + \beta}{2}$};
  \path(3.8,5.2) node()   {\small $\frac{\gamma}{2}$};
  \path(4.3,5) node() 	  {\small $\frac{\gamma}{2}$};
  \path(9.75,-2.1) node() {\small $\frac{\gamma}{2}$};
\end{tikzpicture}
\end{center}

In the following, we will show the derivation of Molweide's equation (2). Although the proof of the other equation (1) can be presented without words~\citep{Kleine88, Wu01}, the proof of Equation (2) needs some forms of wording. To start with, consider $\triangle ABC$ in the figure above. Take an angle bisector through $C$ that cuts $\gamma$ in half. Construct $\triangle ABF$ such that $CD \perp AF$ and $BF \perp AF$. Hence, $\triangle ACE$ is an isosceles triangle, where $E$ is the intersection of sides $BC$ and $AF$. Consequently, $2 \, \measuredangle \, CAD + \gamma = \pi = \alpha + \beta + \gamma$ and thus $\measuredangle \, CAD = (\alpha + \beta)/2 = \measuredangle \, CED$. Furthermore, $\measuredangle \, BAD + \measuredangle \, CAD = \alpha$, so $\measuredangle \, BAD = (\alpha - \beta)/2$. Consider the right triangle $\triangle BEF$. We have $\measuredangle \, BEF = \measuredangle \, CED$ since they are opposed by vertex $E$. So, $\measuredangle \, EBF = \pi/2 - (\alpha + \beta)/2 = \gamma/2$.

Consider $\triangle ABF$, we know that
\begin{equation}
  \cos \frac{1}{2}(\alpha - \beta) = \frac{AF}{AB} = \frac{AF}{c}. \label{3}
\end{equation}
Moreover, consider $\triangle ACD$, $\triangle DCE$ and $\triangle EBF$, respectively. We now have the following relationships:
\begin{align*}
  \sin \frac{1}{2} \gamma &= \frac{AD}{AC} = \frac{AD}{b} 		 &\Longrightarrow&  &AD& = b \sin \frac{1}{2} \gamma \\
						  &= \frac{DE}{CE} = \frac{DE}{b}  		 &\Longrightarrow&  &DE& = b \sin \frac{1}{2} \gamma \\
						  &= \frac{EF}{BE} = \frac{EF}{a - b}    &\Longrightarrow&  &EF& = (a - b) \sin \frac{1}{2} \gamma.
\end{align*}
Furthermore, $AF = AD + DE + EF = (a + b) \sin \frac{1}{2} \gamma$. Substituting this into~\eqref{3}, we obtain the desired Mollweide's formula~\eqref{2}.

As mentioned earlier, Mollweide's formula exhibits cyclical nature. Thus, we can also express Mollweide's formulas~\eqref{1} and~\eqref{2} as follows:
\begin{eqnarray}
  \frac{\sin \frac{1}{2}(\alpha - \beta)}{\cos \frac{1}{2}\gamma} = \frac{a - b}{c} &\qquad&
  \frac{\sin \frac{1}{2}(\beta - \gamma)}{\cos \frac{1}{2}\alpha} = \frac{b - c}{a} \qquad \quad
  \frac{\sin \frac{1}{2}(\gamma - \alpha)}{\cos \frac{1}{2}\beta} = \frac{c - a}{b} \\
  \frac{\cos \frac{1}{2}(\alpha - \beta)}{\sin \frac{1}{2}\gamma} = \frac{a + b}{c} &\qquad&
  \frac{\cos \frac{1}{2}(\beta - \gamma)}{\sin \frac{1}{2}\alpha} = \frac{b + c}{a} \qquad \quad
  \frac{\cos \frac{1}{2}(\gamma - \alpha)}{\sin \frac{1}{2}\beta} = \frac{c + a}{b}.
\end{eqnarray}

One classical reference work on Trigonometry states, `It is absolutely necessary that the computer should know that its results are correct. For this reason all work must be checked'~\citep{Lennes28}. Since Mollweide's formula is useful for checking the result after one solves an oblique triangle, we would like to propose to include this topic when teaching Trigonometry to our students. The equation is particularly meaningful after discussing the Law of Sines and the Law of Cosines for an oblique triangle. The following examples illustrate the usefulness of Mollweide's formula.\\

\noindent {\slshape Example 1.} Use the Law of Sines to solve $\triangle ABC$, where $\alpha = 72^{\circ}$, $\beta = 40^{\circ}$ and $c = 15$. Then use Mollweide's formula to check the answers.

The third angle is $\gamma = 68\,^{\circ}$ and we readily apply the Law of Sines:
\begin{equation}
  \frac{\sin \alpha}{a} = \frac{\sin \beta}{b} = \frac{\sin \gamma}{c}
\end{equation}
which solves $a = 15.39$ and $b = 10.4$. Checking these results using Mollweide's formula~\eqref{1} yields
\begin{align*}
  \textmd{LHS} = \frac{\sin \frac{1}{2}(\alpha - \beta)}{\cos \frac{1}{2}\gamma} = \frac{\sin 16\,^{\circ}}{\cos 34\,^{\circ}} &= 0.3324787260 \qquad \textmd{and} \\
  \textmd{RHS} = \frac{a - b}{c} = \frac{15.37 - 10.4}{15} &= 0.3324787263.
\end{align*}
From this result, we observe that both sides of Mollweide's formula~\eqref{1} have identical values up to $10^{-9}$ error due to computational truncation.\\

\noindent {\slshape Example 2.} Use the Law of Cosines to solve $\triangle ABC$, where $a = 10$, $b = 4$ and $\gamma = 40^{\circ}$. Use Mollweide's formula to check the answers.

One of the Law of Cosines reads as follows:
\begin{equation}
  c^{2} = a^2 + b^{2} - 2ab \cos \gamma = 100 + 16  - 80 \cos 40^{\circ} = 42.72.
\end{equation}
Thus, $c = 7.4$. The other two angles can be found using either by the Law of Sines or the Law of Cosines. We obtain $\alpha = 119.66^{\circ}$ and $\beta = 20.34^{\circ}$. Note that $\alpha$ is an obtuse angle, although our calculator might give an acute angle for it. Checking the answers using Mollweide's formula~\eqref{2} yields
\begin{align*}
  \textmd{LHS} = \frac{\cos \frac{1}{2}(\alpha - \beta)}{\sin \frac{1}{2}\gamma} = \frac{\cos 49.66^{\circ}}{\sin 20^{\circ}} &= 1.892644733 \\
  \textmd{RHS} = \frac{a + b}{c} = \frac{10 + 4}{7.4} &= 1.892644735.
\end{align*}
From this result, we observe that both sides of Mollweide's formula~\eqref{2} have identical values up to $10^{-8}$ error due to computational truncation.

It has been remarked in the Introduction that we can derive the Law of Tangents from Mollweide's formula. It is particularly interesting since the derivation is almost straightforward, in comparison to other derivation procedures. To derive the Law of Tangents, simply divide the two versions of Mollweide's formulas~\eqref{1}--\eqref{2}. We obtain
\begin{equation*}
  \frac{\tan \frac{1}{2}(\alpha - \beta)}{\cot \frac{1}{2} \gamma} = \frac{a - b}{a + b}.
\end{equation*}
Since $\gamma = \pi - (\alpha + \beta)$, then $\cot \frac{1}{2} \gamma = \tan \frac{1}{2}(\alpha + \beta)$ and thus the Law of Tangents is readily obtained.

\section{Conclusion and remark}

In this article, we have revisited Mollweide's formula in the context of teaching Trigonometry to the students at the Foundation (pre-undergraduate) level at The University of Nottingham Malaysia Campus. Although we may tend to skip this material when discussing the law of trigonometric functions, it turns out that this equation is very useful for checking the result after solving oblique triangles and its strong relationship with the Law of Tangents. Depending on the allotted time, it would be interesting and beneficial to investigate Mollweide's formula. This topic could also be regarded as an `enrichment' section, particularly for the students who are interested to dig deeper on the topic. As a remark, Mollweide's formula has found a unique place in Trigonometry as a computational tool. However, it loses its glory as new technology emerges. Nevertheless, we should remember those mathematicians who paved the road ahead of us because their achievements demonstrate the ingenuity of the human mind~\citep{Wu07}.

\section*{\large Acknowledgement}
The author would like to acknowledge Dr Hoo Ling Ping, Dr Mohamad Rafi bin Segi Rahmat and Balrama Applanaidu (The University of Nottingham Malaysia Campus), Dr Joseph Yeo Kai Kow (National Institute of Education, Singapore) and the anonymous referees for the improvement of this paper. This research is supported by the New Researchers Fund NRF 5035-A2RL20 from The University of Nottingham, University Park Campus, UK and Malaysia Campus, Malaysia.

\bigskip \bigskip

\noindent {\bf Natanael Karjanto}, is an instructor at the Department of Mathematics, University College, Sungkyunkwan University, Natural Science Campus, South Korea. Previously, he was an instructor in the Department of Applied Mathematics at The University of Nottingham Malaysia Campus, Malaysia.
\end{document}